\documentclass[letterpaper,11pt]{article}
\usepackage[margin=1in]{geometry}
\usepackage[bookmarks, colorlinks=true, plainpages = false, citecolor = blue,linkcolor=red,urlcolor = blue, filecolor = blue]{hyperref}
\usepackage{url}
\usepackage{amsmath,amsfonts,amsthm,amssymb,bm, verbatim,dsfont,mathtools}
\usepackage{algorithm,algorithmic,color,graphicx,appendix}
\usepackage{subfigure}
\usepackage{etoolbox}
\usepackage{tikz}
\usepackage{xr,xspace}
\usepackage{todonotes}
\usepackage{paralist}
\usepackage{caption}
\theoremstyle{plain}
\newtheorem{theorem}{Theorem}
\newtheorem{lemma}{Lemma}
\newtheorem{proposition}{Proposition}

\theoremstyle{definition}
\newtheorem{definition}{Definition}

\newtheorem*{remark*}{Remark}
\newtheorem*{example}{Example}

\newcommand{\LL}{\mathcal{L}}

\newcommand\1[1]{\mathbb{I}_{\left\{#1\right\}}}
\newtheorem{objective}{Objective}
\newtheorem{assumption}{Assumption}

\newcommand{\SBM}{{\sf SBM}\xspace}
\newcommand{\GSBM}{{\sf GSBM}\xspace}
\newcommand{\calX}{{\mathcal{X}}}
\newcommand{\calP}{{\mathcal{P}}}
\newcommand{\calL}{{\mathcal{L}}}
\newcommand{\eg}{e.g.\xspace}

\begin{document}

\title{Edge Label Inference in Generalized Stochastic Block Models: \\ from Spectral Theory to Impossibility Results}

\date{\today}

\author{Jiaming Xu\and Laurent Massouli\'{e} \and Marc Lelarge\thanks{J. X. is with
the Department of ECE, University of Illinois at Urbana-Champaign, \texttt{jxu18@illinois.edu}. L. M. is with Microsoft Research-Inria Joint Centre, \texttt{laurent.massoulie@inria.fr}. M. L. is with INRIA-ENS, \texttt{marc.lelarge@ens.fr}  }}

\maketitle
\begin{abstract}
The classical setting of community detection consists of networks exhibiting a clustered structure.
To more accurately model real systems we consider a class of networks (i) whose edges may carry labels and (ii) which may lack a clustered structure. Specifically we assume that nodes possess latent attributes drawn from a general compact space and edges between two nodes are randomly generated and labeled according to some unknown distribution as a function of their latent attributes.
Our goal is then to infer the edge label distributions from a partially observed network.
We propose a computationally efficient spectral algorithm and show it allows for asymptotically correct inference when the average node degree could be as low as logarithmic in the total number of nodes. Conversely, if the average node degree is below a specific constant threshold, we show that no algorithm can achieve better inference than guessing without using the observations. As a byproduct of our analysis, we show that our model provides a general procedure to construct random graph models with a spectrum asymptotic to a pre-specified eigenvalue distribution such as a power-law distribution.
\end{abstract}

\section{Introduction}
Detecting communities in networks has received a large amount of attention and has found numerous applications across various disciplines including
physics, sociology, biology, statistics, computer science, etc (see the exposition \cite{Fortunato10} and the references therein). Most previous work assumes networks can be divided into groups of nodes with dense connections internally and sparser connections between groups, and considers random graph models with some underlying cluster structure such as the stochastic blockmodel (\SBM), a.k.a.\ the {\it planted partition model}. In its simplest form, nodes are partitioned into clusters, and any two nodes are connected by an edge independently at random with probability $p$ if they are in the same cluster and with probability $q$ otherwise. The problem of cluster recovery under the \SBM has been extensively studied and many efficient algorithms with provable performance guarantees have been developed (see \eg, \cite{ChenXu14} and the references therein).

Real networks, however, may not display a clustered structure; the goal of community detection should then be redefined.
As observed in \cite{Heimlicher12}, interactions in many real networks can be of various types and prediction of unknown interaction types may have practical merit such as prediction of missing ratings in recommender systems. Therefore an intriguing question arises: Can we accurately predict the unknown interaction types in the absence of a clustered structure?
To answer it, we generalize the \SBM by relaxing the cluster assumption and allowing edges to carry labels. In particular, each node has a latent attribute coming from a general compact space and for any two nodes, an edge is first drawn and then labeled according to some unknown distribution as a function of their latent attributes. Given a partial observation of the labeled graph generated as above, we aim to infer the edge label distributions, which is relevant in many scenarios such as:
\begin{itemize}
\item Collaborative filtering: A recommender system can be represented as a labeled bipartite graph where if a user rates a movie, then there is a labeled edge between them with the label being the rating. One would like to predict the missing ratings based on the observation of a few ratings.
\item Link type prediction: A social network can be viewed as a  labeled graph where if a person knows another person, then there is a labeled edge between them with the label being their relationship type (either friend or colleague). One would like to predict the unknown link types based on the few known link types.
\item Prediction of gene expression levels:  A DNA microarray can be looked as a a labeled bipartite graph where if a gene is expressed in a sample, then there is a labeled edge between them with the label being the expression level. One would like to predict the unobserved expression level based on the few observed expression levels.
\end{itemize}

\subsection{Problem formulation}\label{Section:Model}
The generalized stochastic blockmodel (\GSBM) is formally defined by seven parameters $n$, $\calX,$ $P,$ $B,$ $\calL,$ $\mu,$ $\omega$, where $n$ is a positive integer; $\cal X$ is a compact space endowed with the probability measure $P$; $B: \calX \times \calX \to [0,1]$ is a function symmetric in its two arguments; $\calL$ is a finite set with $\calP(\calL)$ denoting the set of probability measures on it; $\mu: \calX \times \cal X \to \calP(\calL)$ is a measure-valued function symmetric in its two arguments; $\omega$ is a positive real number.

\begin{definition}
Suppose that there are $n$ nodes indexed by $i \in
\{1,\ldots,n\}$. Each node $i$ has an attribute $\sigma_i$ drawn in an
i.i.d.\ manner from the distribution $P$ on ${\mathcal{X}}$.
A random labeled graph  is generated based on $\sigma$: For each pair of nodes $i,j$, independently of all others, we draw an edge between them with probability $B_{\sigma_i, \sigma_j}$; then for each edge $(i,j)$, independently of all others, we label it by $\ell \in \calL$ with probability $\mu_{\sigma_i,\sigma_j}(\ell)$; finally each labeled edge is retained with probability $\omega/n$ and erased otherwise.
\end{definition}
Given a random labeled graph $G$ generated as above, our goal is to infer the edge label distribution $\mu_{\sigma_i,\sigma_j}$ for any pair of nodes $i$ and $j$. To ensure the inference is feasible,  we shall make the following {\em identifiability} assumption: Let $\nu_{x,y}:=B_{x,y} \mu_{x,y}$ and
\begin{equation}\label{identify}
\forall x\neq x'\in {\mathcal{X}},\; \sum_{\ell\in\LL}\int_{\mathcal{X}}|\nu_{x,y}(\ell)-\nu_{x',y}(\ell)|P(dy)>0;
\end{equation}
otherwise $x, x^\prime$ are statistically indistinguishable and can be combined as a single element in $\mathcal{X}$.
We emphasize that the model parameters $(\calX, P, B, \calL, \mu)$ are
all fixed and do not scale with $n$, while $\omega$ could scale with
$n$.  Notice that $n\omega$ characterizes the total number of observed edge labels and thus can be seen as a measure of ``signal strength''.

\subsection{Main results}
We show that it is possible to make meaningful inference of edge label distributions without knowledge of any model parameters in the relatively ``sparse'' graph regime with $\omega=\Omega(\log n)$. In particular, we propose a computationally efficient spectral algorithm with a random weighing strategy. The random weighing strategy assigns a random weight to each label and constructs a weighted adjacency matrix of the label graph.
The spectral algorithm embeds the nodes into a finite, low dimensional Euclidean space based on the leading eigenvectors of the weighted adjacency matrix and uses the empirical frequency of labels on the local neighborhood in the Euclidean space to estimate the underlying true label distribution.

In the very ``sparse'' graph regime with $\omega=O(1)$, since there
exist at least $\Theta(n)$ isolated nodes without neighbors and to
infer the edge label distribution between two isolated nodes the
observed labeled graph $G$ does not provide any useful information, it
is impossible to make meaningful inference for at least a positive
fraction of node pairs. Moreover, we show that it is impossible to
make meaningful inference for any randomly chosen pair of nodes
when $\omega$ is below a specific non-trivial threshold.

As a byproduct of our analysis, we show how the GSBM can generate random graph models with a spectrum asymptotic to a pre-specified eigenvalue distribution such as e.g.\ a power law by appropriately choosing model parameters based on some Fourier analysis.

\subsection{Related work}
Below we point out some connections of our model and results to prior work. More detailed comparisons are provided after we present the main theorems.

\paragraph*{The \SBM and spectral methods} If the node attribute space $\calX$ is a finite set and no edge label is available, then the \GSBM reduces to the classical \SBM with finite number of blocks. The spectral method and its variants are widely used to recover the underlying clusters under the \SBM, see, \eg, \cite{McSherry01,Coja-oghlan10,mastom11,Yu11,ChaudhuriGT12}. However, the previous analysis relies on the low-rank structure of the edge probability matrix. In contrast, the edge probability matrix under the \GSBM is not low-rank, and our analysis is based on establishing a correspondence between the spectrum of a compact operator and the spectrum of a weighted adjacency matrix (see Proposition \ref{PropSpectrum}). Similar connection appears before in the context of data clustering considered in \cite{vonluxburg-bousquet-belkin}, where a graph is constructed based on observed attributes of nodes and clustering based on the graph Laplacian is analyzed. In contrast our setup does not assume the observation of node attributes. Also in our case the observed graphs could be very sparse,
while the graphs considered in~\cite{vonluxburg-bousquet-belkin} are dense. 

\paragraph*{Latent space model} If the node attribute space $\calX$ is a finite-dimensional Euclidean space and no edge label is present, then the \GSBM reduces to the latent space model, proposed in (\cite{Handcock02latentspace,Handcock07}).
 If we further assume the node attribute space $\mathcal{X}$ is the probability simplex endowed with Dirichlet distribution with a parameter $\alpha$, and $B$ is a bilinear function, then the \SBM reduces to the mixed membership \SBM proposed in \cite{airoldi2008mixed}, which is a popular model for studying the overlapping community detection problem.

\paragraph*{Exchangeable random graphs}
If we ignore the edge labels, the \GSBM fits
exactly into the framework of ``exchangeable random graphs'' and the edge probability function
$B$ is known as ``graphon'' (see \eg, \cite{Airoldi13} and the references therein).
It is pointed out in \cite{Bicke09} that some known functions can be used to approximate the graphon, but
no analysis is presented. Our spectral algorithm approximates the graphon using the eigenfunctions and
The exchangeable random graph models with constant average node degrees has been studied in \cite{Bollobas07}, but the focus there is on the phase transition for the emergence of the giant connected component.

\paragraph*{Phase transition if $\omega=O(1)$}
There is an emerging line of works \cite{Decelle11,Mossel12,Mossel13,Massoulie13,Heimlicher12,Jiaming13} that try to identify the sharp phase transition threshold for positively correlated clustering in the regime with a bounded average node degree.
All previous rigorous results focus on the two communities case, while \cite{Decelle11} gives detailed predictions about the phase transition thresholds in the more general case with multiple communities. Here with multiple communities we identify a threshold below which positively correlated clustering is impossible. However, our threshold is not sharp.

\subsection{Notation}
For two discrete probability distributions $\mu$ and $\nu$ on $\mathcal{L}$, let $\|\mu -\nu\|_{\rm TV}:=\frac{1}{2}\sum_{\ell \in \calL} | \mu(\ell) - \nu(\ell) |$ denote the total variation distance. Throughout the paper, we say an event occurs ``a.a.s.'' or ``asymptotically almost surely'' when it occurs with a probability tending to one as $n \to \infty$. We use the standard big O notation. For instance, for two sequences $\{a_n\},\{b_n\}$, $a_n \sim b_n$ means $\lim_{n \to \infty} \frac{a_n}{b_n}=1$.

\section{Spectral reconstruction if $\omega=\Omega(\log n)$} \label{Section:Dense}

Let $A \in \{0,1\}^{n \times n}$ denote the adjacency matrix of $G$ and $L_{ij} \in \LL$ denote the label of edge $(i,j)$ in $G$. Our goal reduces to infer $\mu_{\sigma_i,\sigma_j}$ based on $A$ and $L$. In this section, we study a polynomial-time algorithm  based on the spectrum of a suitably weighted adjacency matrix.
The detailed description is given in Algorithm \ref{alg:spectral} with four steps.

Step 1 defines the weighted adjacency matrix $\tilde{A}$ using a random weighing function $W$ of edge labels. Step 2 extracts the top $r$ eigenvalues and eigenvectors of $\tilde{A}$
for a given integer $r$. Step 3 embeds $n$ nodes in $\mathbb{R}^r$ based on the spectrum of $\tilde{A}$. Step 4 constructs an estimator of $\mu_{\sigma_i,\sigma_j}$ using the empirical label distribution on the edges between node $j$ and nodes in the local neighborhood of node $i$.
Note that the random weight function $W$ chosen in Step 1 is the key to exploit the labeling information encoded in $G$. If $\nu_{x,y}$ were known, better deterministic weight function could be chosen to allow for sharper estimation, e.g.\ (\cite{Jiaming13}). However, no a priori deterministic weight function could ensure consistent estimation irrespective of $\nu_{x,y}$. The function $h_{\epsilon} (x):=\min(1,\max(0,2-x/\epsilon))$ used in Step 4 is a continuous approximation of the indicator function $\1{x \le \epsilon}$ such that $h_{\epsilon}(x)=1$ if $x \le \epsilon$ and $h_{\epsilon}=0$ if $x\ge 2 \epsilon$.

\begin{algorithm}
\caption{Spectral Algorithm ($A,L,r,\epsilon$)}\label{alg:spectral}
\begin{algorithmic}[1]
\STATE (Random Weighing) Let $W:\LL \to [0,1]$ be a random weighing function, with i.i.d.\ weights $W(\ell)$ uniformly distributed on $[0,1]$. Define the weighted adjacency matrix as  $\tilde{A}_{ij} = W(L_{ij})$ if $A_{ij}=1$; otherwise $\tilde{A}_{ij}=0$.

\STATE (Spectral Decomposition) For a given positive integer $r$, extract the $r$ largest eigenvalues of $\tilde{A}$ sorted in decreasing order $|\lambda_1^{(n)}|\ge |\lambda_2^{(n)}| \ge \cdots \ge |\lambda_r^{(n)}| $ and the corresponding eigenvectors with unit norm $v_1,v_2,\ldots,v_r\in\mathbb{R}^n$.

\STATE (Spectral Embedding) Embed the $n$ nodes in $\mathbb{R}^r$ by letting
\begin{align}
z_i:=\sqrt{n}\left(\frac{\lambda_1^{(n)}}{\lambda_1^{(n)}}v_1(i),\ldots,\frac{\lambda_r^{(n)}}{\lambda_1^{(n)}}v_r(i)\right). \label{eq:SpectralEmbedding}
\end{align}
\STATE (Label Estimation) For a given small positive parameter $\epsilon$, define the estimator $\hat{\mu}_{ij}$ of $\mu_{\sigma_i,\sigma_j}$ by letting
\begin{equation}\label{estimator}
\hat{\mu}_{ij}(\ell):=\frac{\sum_{i'} h_\epsilon(||z_{i'}-z_i||_2) \1{L_{i'j}=\ell}}{\epsilon+\sum_{i'} h_\epsilon( ||z_{i'}-z_i||_2) A_{i'j}}.
\end{equation}
Define the estimator $\hat{B}_{ij}$ of $\frac{\omega}{n} B_{\sigma_i,\sigma_j}$ by letting
\begin{align}\label{eq:edgeprobabilityestimator}
\hat{B}_{ij}(\ell):=\frac{\sum_{i'} h_\epsilon(||z_{i'}-z_i||_2) A_{i'j} }{\epsilon+ \sum_{i'} h_\epsilon( ||z_{i'}-z_i||_2)}.
\end{align}
\end{algorithmic}
\end{algorithm}

Our performance guarantee of Spectral Algorithm \ref{alg:spectral} is stated in terms of the spectrum of the integral operator defined as
\begin{equation}\label{operator}
Tf(x):=\int_{\mathcal{X}}K(x,y)f(y) P(dy),
\end{equation}
where the symmetric kernel $K$ is defined by
\begin{equation}
K(x,y):=\sum_{\ell} W(\ell)\nu_{x,y}(\ell) \in [0,|\LL|].
\end{equation}
Since $K$ is bounded, the operator $T$, acting on the function space $L_2(P)$, is {\em compact} and therefore admits a discrete spectrum with finite multiplicity of all of its non-zero eigenvalues (see e.g.\ \cite{kato66} and \cite{vonluxburg-bousquet-belkin}). Moreover, any of its eigenfunctions is continuous on $\mathcal{X}$. Denote the eigenvalues of operator $T$ sorted in decreasing order by $|\lambda_1| \ge |\lambda_2| \ge \cdots$ and its corresponding eigenfunctions with unit norm by $\phi_1, \phi_2, \cdots$. Define
\begin{equation}\label{identify2}
d^2(x,x'):=\;\int_{\mathcal{X}} |K(x,y)-K(x',y)|^2 P(dy).
\end{equation}
It is easy to check that with probability 1 with respect to the random choices of $W(\ell)$,
by the identifiability condition~\eqref{identify}, $d(x,x')>0$ for all
$x \neq x' \in \mathcal{X}$. By Minkowski inequality, $d(x,x')$
satisfies the triangle inequality. Therefore, $d(x,x')$ is a distance
on $\mathcal{X}$. By the definition of $\lambda_k$ and $\phi_k$, we
have (the following serie converges in $L_2(P\times P)$, see Chapter
V.4 in \cite{kato66}):
\begin{align}
K(x,y)=\sum_{k=1}^\infty \lambda_k \phi_k(x) \phi_k(y), \label{eq:kernalapprox}
\end{align}
and thus $d^2(x,x')=\sum_{k=1}^\infty \lambda_k^2 \left(\phi_k(x)-\phi_k(x')\right)^2$.

To derive the performance guarantee of Spectral Algorithm \ref{alg:spectral}, we make the following continuity assumption on $\nu_{x,y}$. Similar continuity assumptions appeared before in the literature on the latent space model and the exchangeable random graph model (see \eg, \cite[Section 4.4]{Chattergee12} and \cite[Section 2.1]{Airoldi13}).
\begin{assumption}\label{assumptioncontinuous}
For every $\ell\in\LL$, $\nu_{x,y}(\ell)$ is continuous on $\mathcal{X}^2$, hence by compactness of ${\mathcal{X}}$ uniformly continuous. Let $\psi(\cdot)$ denote a modulus of continuity of all functions $(x,y)\to\nu_{x,y}(\ell)$ and $(x,y)\to B_{x,y}$. That is to say, for all $x,x',y,y'$,
$$
|B_{x,y}-B_{x',y'}|\le \psi(d(x,x')+d(y,y'))
$$
and similarly for $\nu_{x,y}(\ell)$.
\end{assumption}

Let $\epsilon_r:=\sum_{k>r}\lambda_k^2$ for a fixed integer $r$, characterizing the tail of the spectrum of the operator $T$. The following theorem gives an upper bound of the estimation error of $\hat{\mu}_{ij}$ for most pairs $(i,j)$ in terms of $\epsilon_r$ and $\epsilon$.
\begin{theorem}\label{thm1}
Suppose Assumption~\ref{assumptioncontinuous} holds. Assume that $ \omega \ge C \log n $ for some universal positive constant $C$ and $r$ chosen in Spectral Algorithm \ref{alg:spectral} satisfies $|\lambda_r|>|\lambda_{r+1}|$.
Then a.a.s.\ the estimators $\hat{\mu}$ and $\hat{B}$ given in Spectral Algorithm \ref{alg:spectral} satisfy
\begin{align}
B_{\sigma_i, \sigma_j} | \hat{\mu}_{ij} (\ell)-\mu_{\sigma_i, \sigma_j}(\ell) |  & \le  2 \psi(2|\lambda_1|\epsilon) + \frac{1}{|\lambda_1|^2\epsilon^2} \frac{ \sqrt{\epsilon_r}  }{ \int_{\mathcal{X} } h_{|\lambda_1| \epsilon}(d(\sigma_i,x) ) P(dx)}:= \eta, \; \forall \ell \in \calL, \nonumber \\
| \hat{B}_{ij}- \frac{\omega}{n} B_{\sigma_i, \sigma_j} | &  \le \frac{\omega}{n} \eta,  \label{eq:errorbound}
\end{align}
for a fraction of at least $(1-\sqrt{\epsilon_r})$ of all possible pairs  $(i,j)$ of nodes.
\end{theorem}
Note that if $\epsilon_r$ goes to $0$, the second term in $\eta$ given by~\eqref{eq:errorbound} vanishes, and $\eta$ simplifies to $  2\psi(2|\lambda_1|\epsilon)$ which goes to $0$ if $\epsilon$ further goes to $0$.  In the case where $B_{\sigma_i,\sigma_j}$ is strictly positive, Theorem \ref{thm1} implies that the estimation error of the edge label distribution goes to $0$ as successively $\epsilon_r$ and $\epsilon$ converge to $0$.
Note that $\epsilon$ is a free parameter chosen in Spectral Algorithm \ref{alg:spectral} and can be made arbitrarily small if $\epsilon_r$ is sufficiently small.
The parameter $\epsilon_r$ measures how well the compact space $\mathcal{X}$ endowed with measure $P$ can be approximated by $r$ discrete points, or equivalently how well our general model can be approximated by the labeled stochastic block model with $r$ blocks.
The smaller $\epsilon_r$ is, the more structured, or the more ``low-dimensional'' our general model is. In this sense,  Theorem \ref{thm1} establishes an interesting connection between the estimation error and the structure present in our general model.

A key part of the proof of Theorem \ref{thm1} is to show that for any fixed $k$, the normalized $k$-th largest eigenvalue $\lambda_k^{(n)}/\lambda_1^{(n)}$  of the weighted adjacency matrix $\tilde{A}$ asymptotically converges to $\lambda_k/\lambda_1$ where $\lambda_k$ is the $k$-th eigenvalue of integral operator $T$, and this is precisely why our spectral embedding given by~\eqref{eq:SpectralEmbedding} is defined in a normalized fashion.
The following simple example illustrates how we can derive closed form expressions for the spectrum of integral operator $T$.
\begin{example} \label{ex1}
Take ${\mathcal{X}}=[0,1]$ and $P$ as the Lebesgue measure. Assume unlabeled edges. Let $B_{x,y}=g(x-y)$ where $g$ is an even (i.e.\ $g(-\cdot)=g(\cdot)$), 1-periodic function. Denote its Fourier series expansion  by $g(x)=\sum_{k\ge 0}g_k \cos(2\pi k x)$. For instance, if $g(x)=|x|$ for $x\in[-1/2,1/2]$, then $g_0=1/4$ and $g_k=[(-1)^k-1]/(\pi^2 k^2)$ for
$k \ge 1$. If $g(x)=\1{ -1/4 \le x \le 1/4}$ for $ x \in [-1/2,1/2]$, then $g_0=1/2$ and $g_k=2\sin(\pi k/2)/(\pi k)$ for $k \ge 1$.
\end{example}
For the above example, $Tf= g \ast f$ where $\ast$ denotes convolution. Fourier series analysis entails that $\lambda_k$ must coincide with Fourier coefficient $g_0$ or $g_k/2$ for $k\ge1$ ($g_k/2$ appearing twice in the spectrum of $T$). This example thus gives a general recipe for constructing random graph models with spectrum asymptotic to a pre-specified eigenvalue profile. For $g(x)=|x|$ on $[-1/2,1/2]$, we find in particular that $\lambda_1=1/4$ and $|\lambda_{2k}|=|\lambda_{2k+1}|=1/(\pi^2 (2k-1)^2)$, which is a power-law spectrum with the decaying exponent being $2$. For $g(x)=\1{ -1/4 \le x \le 1/4}$ on $[-1/2,1/2]$, $\lambda_1=1/2$ and $|\lambda_{2k}|=|\lambda_{2k+1}|=1/(\pi (2k-1) )$, which is a power-law spectrum with the decaying exponent being $1$.

\paragraph*{Comparisons to previous work}
Theorem \ref{thm1} provides the first theoretical result on inferring edge label distributions to our knowledge.
For estimating edge probabilities,
Theorem \ref{thm1} implies or improves the best known results in several special cases.

For the \SBM with finite $r$ blocks, $\epsilon_r$ is zero. By choosing $\epsilon$
sufficiently small in Theorem \ref{thm1}, we see that our spectral method is
asymptotically correct if $\omega=\Omega(\log n)$, which
matches with best known bounds (see \eg, \cite{ChenXu14} and the references therein). For the mixed membership \SBM with finite $r$ blocks, the best known performance guarantee given by \cite{anandkumar2013tensormixed}
needs $\omega $ to be above the order of several $\log n$ factors,
while Theorem \ref{thm1} only needs $\omega$ to be the order of $\log n$. However,
Theorem \ref{thm1} requires the additional spectral gap assumption
and needs $\epsilon_r$ to vanish. Also, notice that Theorem \ref{thm1} only applies to the setting where the edge probability $p$ within the community
exceeds the edge probability $q$ across two different communities by a constant factor,
while the best known results in \cite{ChenXu14,anandkumar2013tensormixed} apply to the general setting with any $r,p,q$.

For the latent space model, \cite{Chattergee12} proposed a universal singular
value thresholding approach and showed that the edge probabilities can
be consistently estimated if $\omega \ge n^{\frac{k}{k+2}}$ with some Lipschitz
condition on $B$ similar to Assumption~\ref{assumptioncontinuous}, where $k$ is the dimension of the node attribute space. Our
results in Theorem \ref{thm1} do not depend on the dimension of the node attribute space and
only need $\omega$ to be on the order of $\log n$.

For the exchangeable random graph models, a singular value thresholding approach is shown in \cite{Chattergee12} to estimate the graphon consistently. More recently, \cite{Airoldi13} shows that the graphon can be consistently estimated using the empirical frequency of edges in local neighborhoods, which are constructed by thresholding based on the pairwise distances between different rows of the adjacency matrix.
All these previous works assume the edge probabilities are constants.  In contrast, Theorem \ref{thm1} applies to much
sparser graphs with edge probabilities could be as low as $\log n/n$.

\section{Impossibility if $\omega=O(1)$} \label{Section:Sparse}

We have seen in the last section that Spectral Algorithm \ref{alg:spectral} achieves asymptotically correct inference of edge label distributions so long as $\omega= \Omega(\log n)$ and $\epsilon_r =0$. In this section, we focus on the sparse regime where $\omega$ is a constant, i.e., the average node degree is bounded and the number of observed edge labels is only linearly in $n$. We identify a non-trivial threshold under which  it is fundamentally impossible to infer the edge label distributions with an accuracy better than guessing without using the observations.

To derive the impossibility result, let us consider a simple scenario where the compact space $\mathcal{X}=\{1, \ldots, r \}$ is endowed with a uniform measure $P$, $B_{x, y} = \frac{a }{a+b} $ if $x=y$ and $B_{x,y}=\frac{b}{a+b}$ if $x \neq y$ for two positive constants $a, b$, and $\mu_{x,y}=\mu$ if $x=y$ and $\mu_{x,y}=\nu$ if $x \neq y$ for two different discrete probability measures  $\mu, \nu$ on $\mathcal{L}$. Since $\omega$ is a constant, the observed labeled graph $G$ is sparse and has a bounded average degree. Similar to the Erd\H{o}s-R\'enyi random graph, there are at least $\Theta(n)$ isolated nodes without neighbors. To infer the edge label distribution between two isolated nodes, the observed labeled graph $G$ does not provide any useful information and thus it is impossible to achieve the asymptotically correct inference of edge label distribution for two isolated nodes. Hence we resort to a less ambitious objective.

\begin{objective}\label{objective}
Given any two  randomly chosen nodes $i$ and $j$, we would like to correctly determine whether the label distribution is $\mu$ or $\nu$ with probability strictly larger than $1-1/r$, which is achievable by always guessing $\nu$ and is the best one can achieve if no graph information available.
\end{objective}
Note that if Objective \ref{objective} is not achievable, then the expected estimation error is at least $\frac{1}{2r} \| \mu- \nu \|_{\rm TV}$. One might think that we can always achieve Objective \ref{objective} as long as $\omega>0$ such that the graph contains a giant connected component, because the labeled graph $G$ then could provide extra information. It turns out that this is not the case. Define
\begin{align}
\tau= \frac{1}{r(a+b)} \sum_{ \ell \in \mathcal{L} } |a\mu(\ell)-b\nu(\ell)|. \label{DefReconstructionThresholdMultiple}
\end{align}
Let $\omega_0=1/\tau$ and $\omega_c= \frac{r(a+b)}{a+(r-1)b}$. Then by definition of $\tau$, we have $\omega_0 > \omega_c$. Note that when $\omega>\omega_c$, the average node degree is larger than one, and thus similar to Erd\H{o}s-R\'enyi random graph, $G$ contains a giant connected component. The following theorem shows that Objective \ref{objective} is fundamentally impossible if $\omega < \omega_0$ where $\omega_0$ is strictly above the threshold $\omega_c$ for the emergence of the giant connected component.

\begin{theorem} \label{ThmNonReconstructionMultiple}
If $ \omega< \omega_0$, then for any two randomly chosen nodes $\rho$ and $v$,
\begin{align}
\forall x,y \in\{1, \ldots, r\},\; \mathbb{P} (\sigma_\rho=x | G, \sigma_v= y ) \sim \frac{1}{r} \text{ a.a.s }. \nonumber
\end{align}
\end{theorem}
The above theorem implies that it is impossible to correctly determine whether two randomly chosen nodes have the same attribute or not with probability larger than $1-1/r$ and thus Objective \ref{objective} is fundamentally impossible. In case $a \neq b$, it also implies that we cannot correctly determine whether the edge probability between nodes $i$ and $j$ is $\frac{a}{a+b}$ or $\frac{b}{a+b}$ with probability strictly larger than $1-1/r$. This indicates the need for a minimum number of observations in order to exploit the information encoded in the labeled graph.



\paragraph*{Comparisons to previous work}
To our knowledge, Theorem \ref{thm1} provides the first impossibility result on inferring edge label distributions and node attributes in the case with multiple communities. The previous work focuses on the case with two community case. If $r=2$ and no edge label is available, it is conjectured in \cite{Decelle11} and later proved in \cite{Mossel12,Mossel13,Massoulie13} that the positively correlated clustering is feasible if and only if $(a-b)^2>2(a+b)$, or equivalently, $\omega> 1/2\tau^2$. If the edge label is available, it is conjectured in \cite{Heimlicher12} that the the positively correlated clustering is feasible if and only if $\omega>1/\tau'$ with
\begin{align*}
\tau' =  \frac{1}{2(a+b)} \sum_{ \ell \in \mathcal{L} } \frac{ (a\mu(\ell)-b\nu(\ell))^2}{ a\mu(\ell)+b\nu(\ell)} \le \tau.
\end{align*}
It is proved in \cite{Jiaming13} that the positively correlated clustering is infeasible  if $\omega<1/\tau'$. Comparing to the previous works, the threshold $1/\tau$ given by Theorem \ref{ThmNonReconstructionMultiple} is not sharp in the special case with two communities.

\section{Numerical experiments}
In this section, we explore the empirical performance of our Spectral Algorithm \ref{alg:spectral} based on Example \ref{ex1}. In particular, suppose $n=1500$ nodes are uniformly distributed over the space $\mathcal{X}=[0,1]$. Let $B_{x,y}=g(x-y)$ where $g$ is even, $1$-periodic and defined by $g(x)=|x|$ for $x \in [-1/2, 1/2]$.
Assume unlabeled edges first.

\begin{figure}[ht]
\centering
\subfigure{%
\includegraphics[width=2.5in]{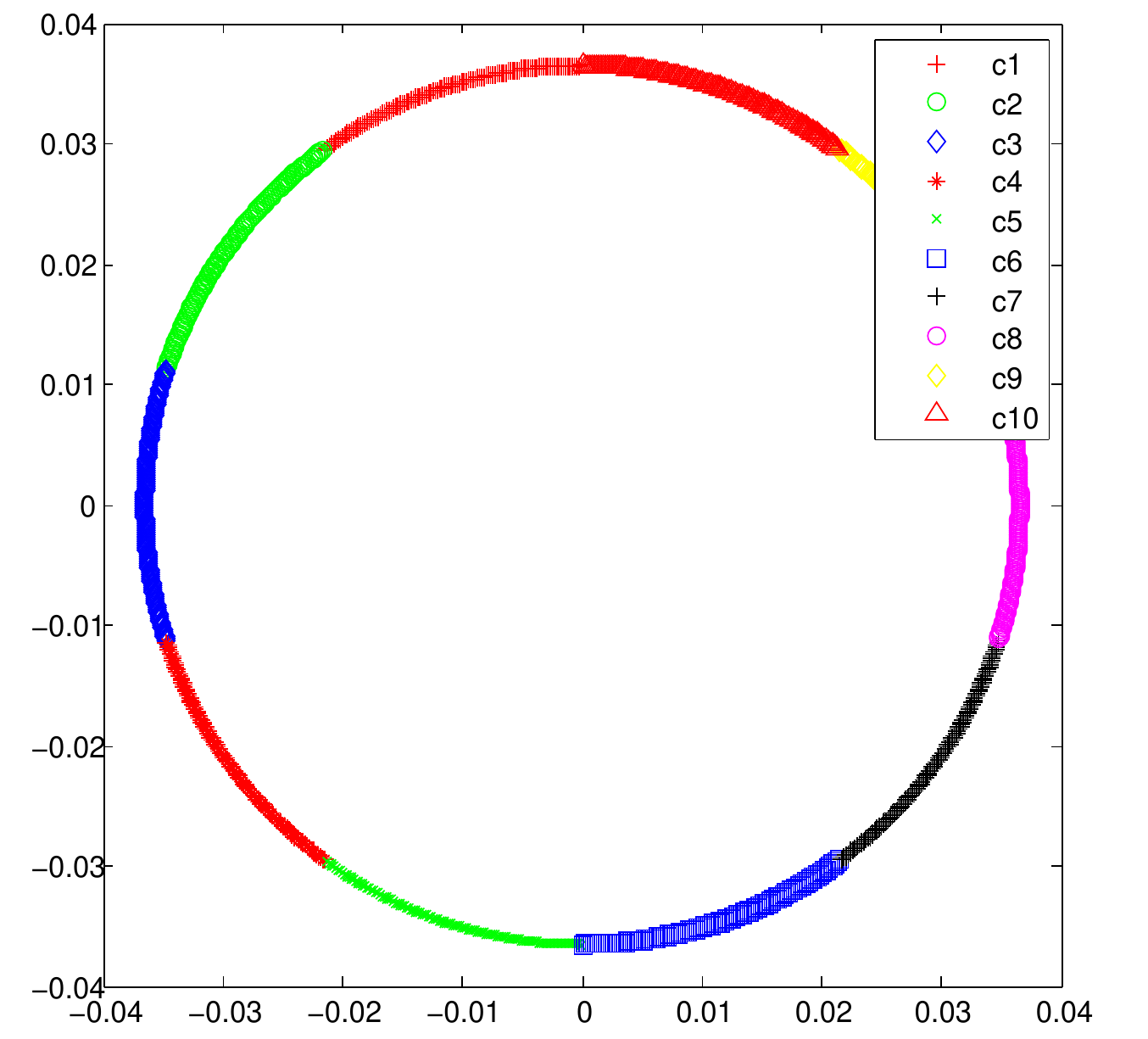}
\label{FigIdeaSpectralEmbedding}}
\quad
\subfigure{
\includegraphics[width=2.5in]{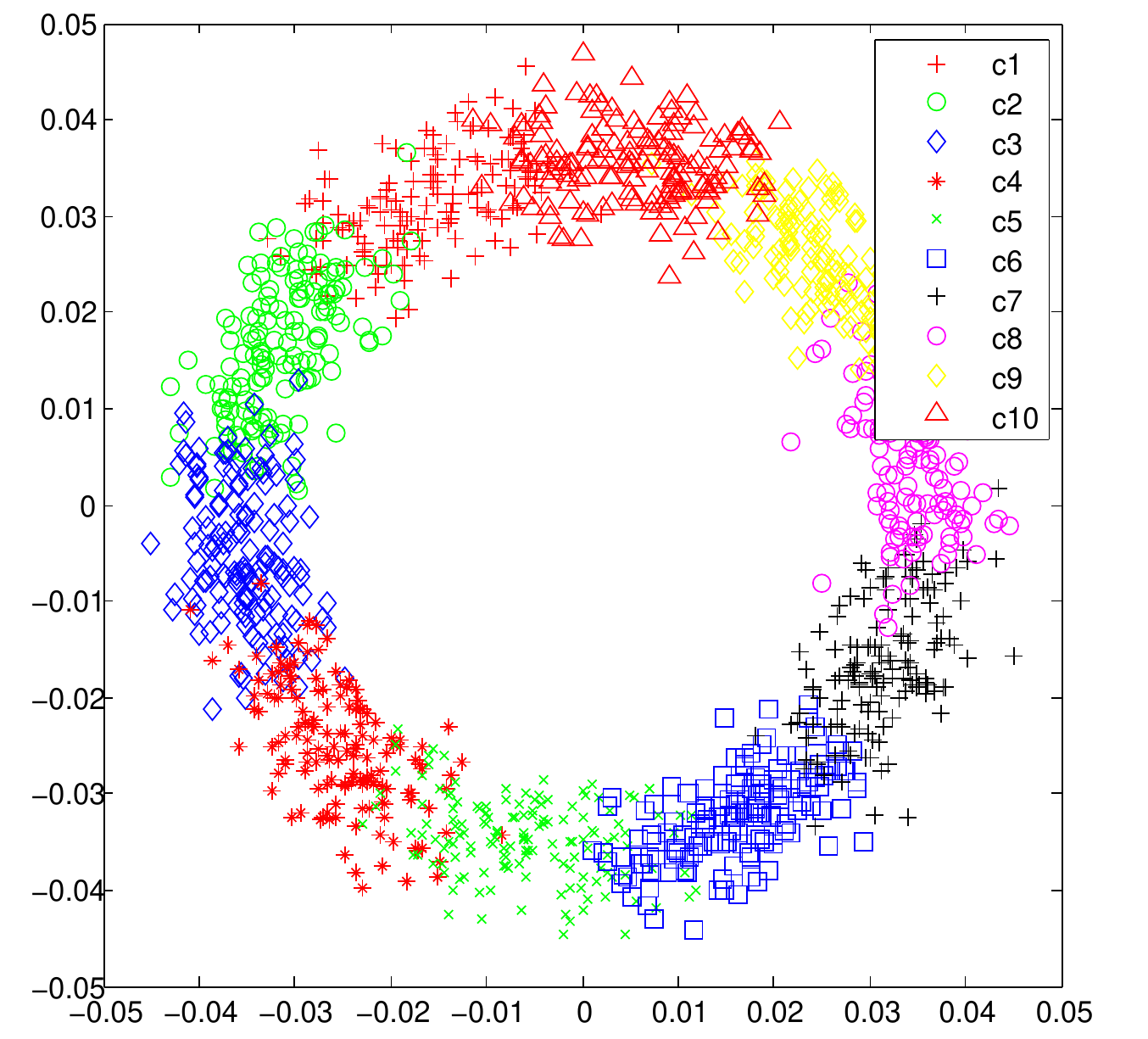}
\label{FigSpectralEmbedding} }
\caption{(a): The spectral embedding given by $f_i$; (b): The spectral embedding given by $z_i$.}
\end{figure}

We simulate the spectral embedding given by Step 3 of Algorithm \ref{alg:spectral} for a fixed observation probability $\omega/n=0.6$.
Pick $r=3$ in Algorithm \ref{alg:spectral}. Note that the eigenvector $v_1$ corresponding to the largest eigenvalue is nearly parallel to the all-one vector and thus does not convey any useful information. Therefore, our spectral embedding of $n$ nodes are based on $v_2$ and $v_3$. In particular, let $z_i=(v_3(i),v_2(i)) \in \mathbb{R}^2$.
As we derived in Section \ref{Section:Dense}, the second and
third largest eigenvalues of operator $T$ are given by $\lambda_2=\lambda_3=-1/\pi^2$, and the corresponding
eigenfunctions are given by $\phi_2(x)= \sqrt{2} \cos(2 \pi x)$ and $\phi_3(x)=\sqrt{2} \sin(2 \pi x)$.
Proposition \ref{PropSpectrum} shows that
$z_i$ asymptotically converges to  $ f_i=\sqrt{\frac{2}{n}} ( \cos( 2 \pi \sigma_i), \sin( 2 \pi \sigma_i) )$.
We plot $f_i$ and $z_i$ in a two-dimensional plane as shown in
Fig.~\ref{FigIdeaSpectralEmbedding} and Fig.~\ref{FigSpectralEmbedding}, respectively.
For better illustration, we divide all nodes into ten groups with different colors, where the $k$-th group consists of nodes with attributes given by $ \frac{1}{n} [ 100(k-1)+1, 100(k-1)+2, \ldots, 100k ]$. As we can see, $z_i$ is close to $f_i$ for most nodes $i$, which coincides with our theoretical finding.

\begin{figure}[ht]
\centering
\subfigure{
\includegraphics[width=2.5in]{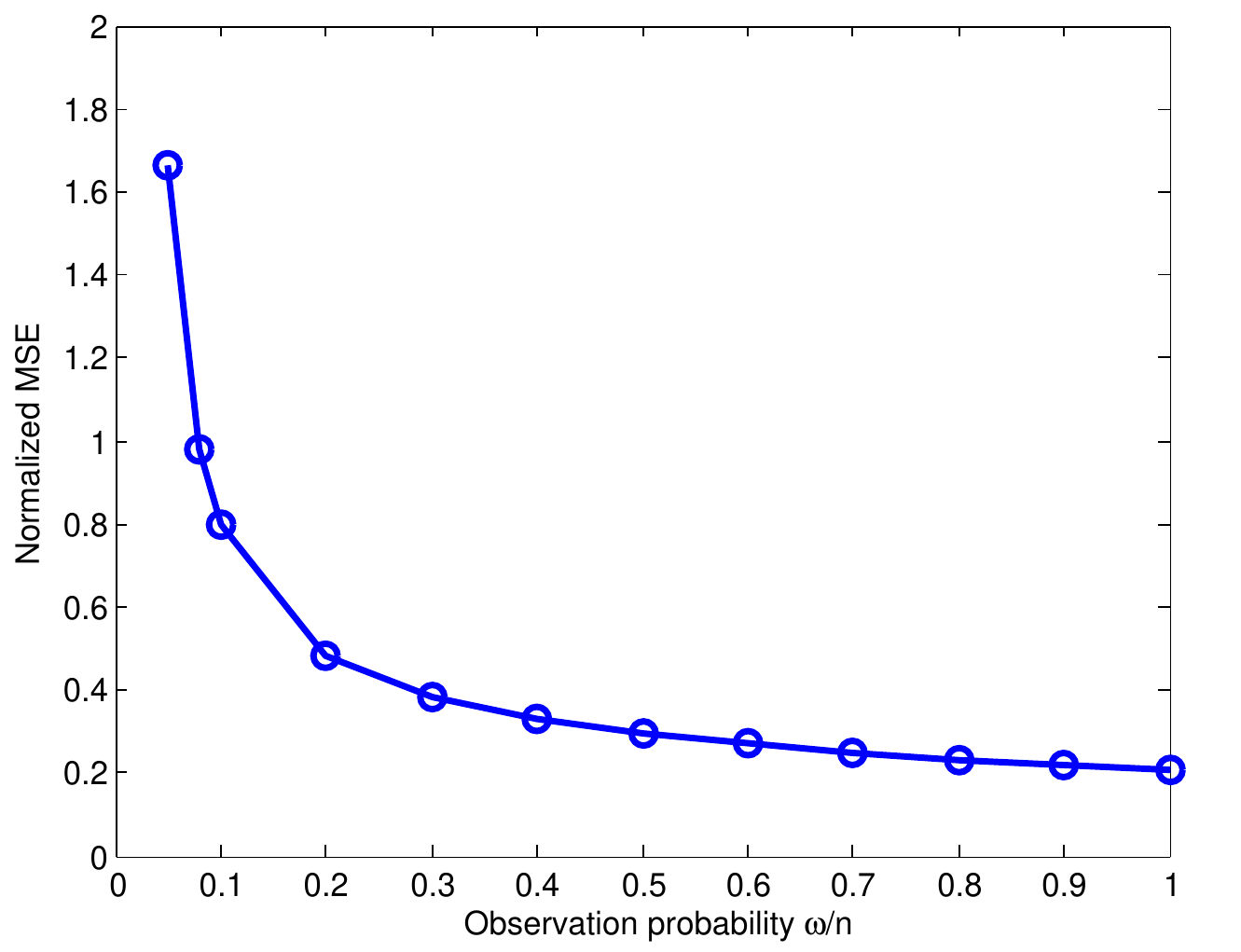}
\label{FigEdgeProbability}}
\quad
\subfigure{%
\includegraphics[width=2.5in]{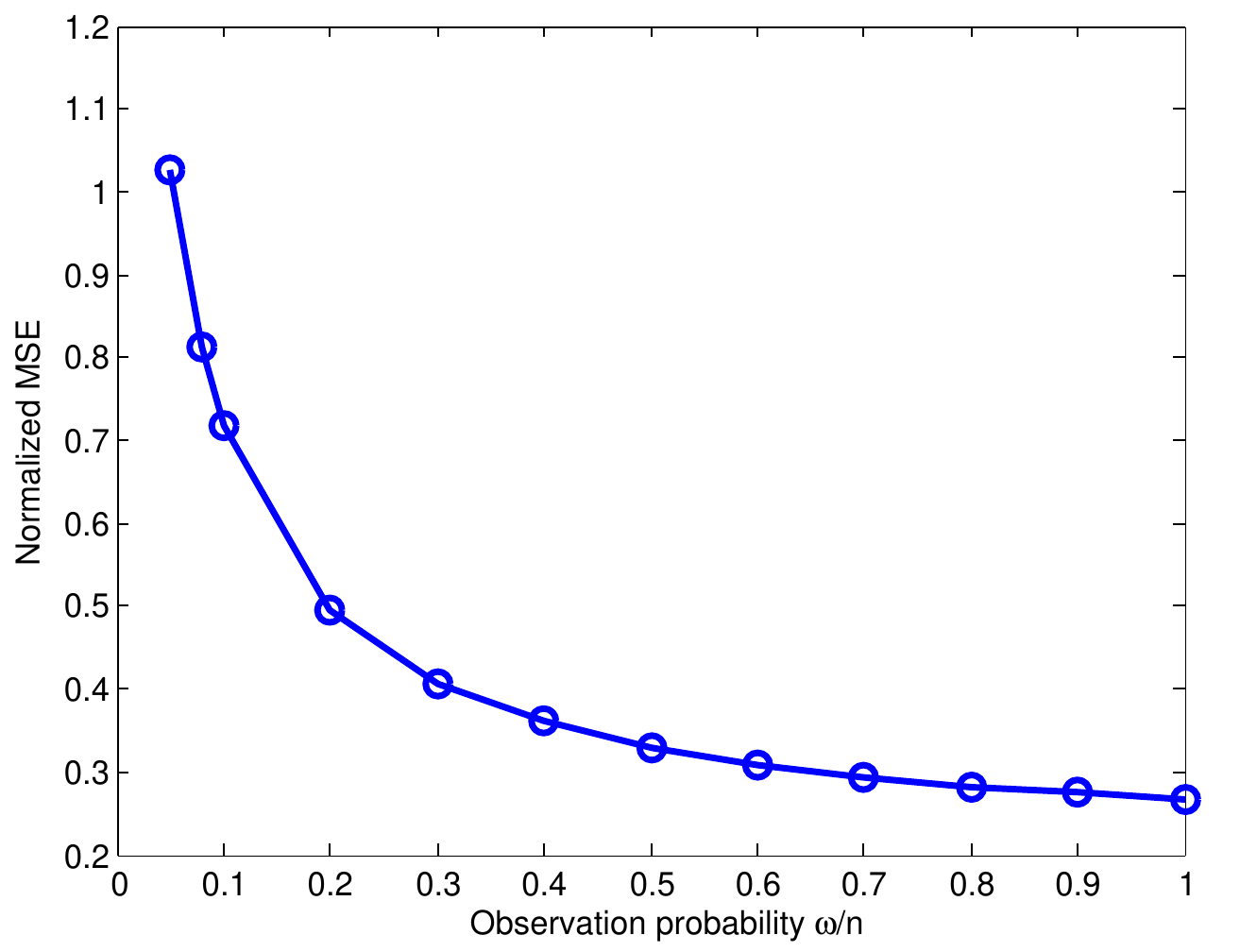}
\label{FigEdgeLabelProbability}
}
\caption{(a): Estimating the observed edge probability; (b):
Estimating the edge label distribution.}
\label{FigEstimation}
\end{figure}
Then we simulate Spectral Algorithm \ref{alg:spectral} on estimating the observed edge probability $\frac{\omega}{n} B_{\sigma_i,\sigma_j}$ between any node pair $(i,j)$ by picking $r=3$ and setting $\epsilon=0.5\text{median} \{ \|z_i-z_j\|\}$.
We measure the estimation error by the normalized mean square error given by $\| \hat{B}- \frac{\omega}{n} B^\ast \|_F /  \| \bar{B}- \frac{\omega}{n} B^\ast  \|_F$,
where $B^\ast$ is the true edge probability defined by $B^\ast_{ij}=B_{\sigma_i,\sigma_j}$; $\hat{B}$ is our estimator defined in~\eqref{eq:edgeprobabilityestimator}; $\bar{B}_{ij}$ is the empirical average edge probability defined by $\bar{B}_{ij}= \sum_{i'} A_{i,i'}/ (n-1).$
Our simulation result is depicted in Fig.~\ref{FigEdgeProbability}.


Next we consider labeled edges with two possible labels $+1$ or $-1$ and $\mu_{x,y}(+1)=2 g(x-y)$. We simulate Spectral Algorithm \ref{alg:spectral} for estimating $\mu_{\sigma_i,\sigma_j}$ between any node pair $(i,j)$ by choosing the weight function as $W(\pm 1)=\pm 1$. We again measure the estimation error by the normalized mean square error given by $\| \hat{\mu} - \frac{\omega}{n} \mu^\ast \|_F /  \| \bar{\mu}- \mu^\ast \|_F$,
where $\mu^\ast$ is the true label distribution defined by $\mu^\ast_{ij} =\mu_{\sigma_i,\sigma_j}$; $\hat{\mu} $  is our estimator defined in~\eqref{estimator}; $\bar{\mu}_{ij}$ is the empirical label distribution defined by $\bar{\mu}_{ij}(+1)= \sum_{i'} \1{L_{ii'}=+1}/ \sum_{i'} A_{ii'}.$
Our simulation result is depicted in Fig.~\ref{FigEdgeLabelProbability}.
As we can see from Fig.~\ref{FigEstimation}, when $\omega/n$ is larger than $0.1$, our spectral algorithm performs better than the estimator based on the empirical average.

\section{Proofs} \label{Proofs}
\subsection{Proof of Theorem \ref{thm1}}


Our proof is divided into three parts. We first establish the asymptotic correspondence
between the spectrum of the weighted adjacency matrix $\tilde{A}$ and the spectrum of the operator $T$ using Proposition \ref{PropSpectrum}.
Then, we prove that the estimator of edge label distribution converges to a limit. Finally, we upper bound the total variation distance between the limit and the true label distribution using Proposition~\ref{Prop2}.

\begin{proposition}\label{PropSpectrum}
Assume that $ \omega \ge C \log n $ for some universal positive constant $C$ and $r$ chosen in Spectral Algorithm \ref{alg:spectral} satisfies $|\lambda_r|>|\lambda_{r+1}|$. Then for $k=1,2,\ldots, r+1$, almost surely $\lambda^{(n)}_k/ \lambda^{(n)}_1 \sim \lambda_k / \lambda_1$. Moreover, for $k=1,2,
\ldots, r$, almost surely there exist choices of orthonormal eigenfunctions $\phi_k$ of operator $T$ associated with $\lambda_k$ such that $\lim_{n\to\infty}\sum_{i=1}^n( v_k(i)-\frac{1}{\sqrt{n}}\phi_k(\sigma_i))^2=0$.
\end{proposition}

By Proposition~\ref{PropSpectrum}, we get the existence of eigenfunctions $\phi_k$ of $T$ associated with $\lambda_k$ such that a.a.s., by letting
$$
f_m:=\left(\frac{\lambda_1}{\lambda_1} \phi_1(\sigma_m),\ldots,\frac{\lambda_r}{\lambda_1} \phi_r(\sigma_m)\right),
$$
we have
\begin{align*}
\sum_{m=1}^n||z_m-f_m ||_2^2= \sum_{m=1}^n \sum_{k=1}^r \left(\sqrt{n}\frac{\lambda^{(n)}_k}{\lambda^{(n)}_1}v_k(m)-\frac{\lambda_k}{\lambda_1}\phi_k(\sigma_m)\right)^2  =o(n).
\end{align*}
By Markov's inequality,
\begin{align*}
\frac{ 1}{n} \left|\{m\in\{1,\ldots,n\}:||z_m-f_m||_2 \ge \delta_n \}\right| \le \frac{\sum_{m=1}^n||z_m-f_m||_2^2}{n \delta_n^2} =\frac{1}{\delta_n^2}o(1).
\end{align*}
Note that $\delta_n$ can be chosen to decay to zero with $n$ sufficiently slowly so that the right-hand side of the above is $o(1)$. We call nodes $m$ satisfying $ \|z_m-f_m\|_2 \ge \delta_n$ ``bad nodes''. Let $\mathcal{I}$ denote the set of ``bad nodes''.
It follows from the last display that $| \mathcal{I}| = o(n)$.  Let $\mathcal{J}$ denote the set of nodes with at least $\gamma_n$ fraction of edges directed towards ``bad nodes'', i.e.,
$$
\mathcal{J}=\{ j : | \{ i \in \mathcal{I} : A_{ij}=1 \}  | \ge \gamma_n | \{ i: A_{ij}=1 \} | \}.
$$
Note that the average node degree in $G$ is $\Theta(\omega)$. Since $\omega \ge C\log n$ by assumption, it follows from the Chernoff bound that the observed node degree is $\Theta(\omega)$ with high probability. Therefore, we can choose $\gamma_n$ decaying to zero while still having $|\mathcal{J} |=o(n)$, i.e., all but a vanishing fraction of nodes have at most $\gamma_n$
fraction of edges directed towards ``bad nodes''.
We have thus performed an embedding of  $n$ nodes in $\mathbb{R}^r$ such that for $m, m' \notin \mathcal{I} $,
\begin{equation}\label{embed0}
||z_m-z_{m'}||_2=\frac{1}{|\lambda_1|}d_r(\sigma_m,\sigma_{m'})+O(\delta_n),
\end{equation}
where pseudo-distance $d_r$ is defined by $
d^2_r(x,x'):=\sum_{k=1}^r\lambda_k^2\left(\phi_k(x)-\phi_k(x')\right)^2.
$

The remainder of the proof exploits this embedding and the
fact that pseudo-distance $d_r$ and distance $d$ are apart by at most $\epsilon_r$ in some suitable sense.
For a randomly selected pair of nodes $(i,j)$, one has a.a.s.\ $i \notin \mathcal{I}$ and $j \notin \mathcal{J}$.  Therefore, node $j$ has
at most $\gamma_n=o(1)$ fraction of edges directed towards ``bad nodes''. Hence, by~\eqref{embed0},
\begin{align}
\sum_{i'} h_{\epsilon} ( ||z_{i'}-z_i||_2) \1{L_{i' j} =\ell }= \sum_{i'} \1{L_{i' j }=\ell}  h_{|\lambda_1|\epsilon} \left( d_r(\sigma_i,\sigma_{i'})+O(\delta_n) \right)+O(\omega \gamma_n), \label{eq_numerator}
\end{align}
and
\begin{align}
\sum_{i'} h_{\epsilon} ( ||z_{i'}-z_i||_2)A_{i'j} =  \sum_{i'}  h_{|\lambda_1|\epsilon}\left( d_r(\sigma_i,\sigma_{i'}) +O(\delta_n) \right) A_{i'j} +O(\omega \gamma_n). \label{eq_denomenator}
\end{align}
The first term in the R.H.S.\ of~\eqref{eq_numerator} is a sum of i.i.d.\ bounded random variables with mean given by
$$
\frac{\omega}{n} \int_{\mathcal{X}} h_{|\lambda_1|\epsilon} \left( d_r(\sigma_i,x)+O(\delta_n) \right) \nu_{x,\sigma_j}(\ell)P(dx).
$$
Since $\omega \ge C \log n$ by assumption, it follows from the Bernstein inequality that a.a.s.\
\begin{align}
\sum_{i'} h_{\epsilon} ( ||z_{i'}-z_i||_2) \1{L_{i' j} = \ell } = &(1+o(1)) \omega \int_{\mathcal{X}} h_{|\lambda_1|\epsilon} \left( d_r(\sigma_i,x)+O(\delta_n) \right) \nu_{x,\sigma_j}(\ell)P(dx)   \nonumber \\
&+ O(\omega \gamma_n). \label{eq_embed1}
\end{align}
The first term in the R.H.S.\ of~\eqref{eq_denomenator} is a sum of i.i.d.\ bounded random variables with mean given by
$$
\frac{\omega}{n} \int_{\mathcal{X}} h_{|\lambda_1|\epsilon}\left( d_r(\sigma_i,x) +O(\delta_n) \right) B_{x,\sigma_j}P(dx).
$$
It again follows from the Bernstein inequality that a.a.s.\
\begin{align}
\sum_{i'} h_{\epsilon} ( ||z_{i'}-z_i||_2)A_{i'j} = & (1+o(1)) \omega \int_{\mathcal{X}} h_{|\lambda_1|\epsilon}\left( d_r(\sigma_i,x) +O(\delta_n) \right) B_{x,\sigma_j} P(dx)  \nonumber  \\
&+ O(\omega \gamma_n). \label{eq_embed2}
\end{align}
Note that $h_\epsilon(x)$ is a continuous function in $x$. Therefore,
$$\lim_{n \to \infty} h_{|\lambda_1|\epsilon}\left( d_r(\sigma_i,x) +O(\delta_n) \right) = h_{|\lambda_1|\epsilon}( d_r(\sigma_i,x) ). $$
By the dominated convergence theorem, it follows from~\eqref{estimator},~\eqref{eq_embed1},~\eqref{eq_embed2} that a.a.s.
\begin{align}
\hat{\mu}_{i,j} (\ell) \sim \frac{   \int_{\mathcal{X} } h_{|\lambda_1| \epsilon}(d_r(\sigma_i,x) ) \nu_{x,\sigma_j} (\ell) P(dx)  }{   \int_{\mathcal{X}} h_{|\lambda_1| \epsilon}(d_r(\sigma_i,x) ) B_{x,\sigma_j} P(dx) }:=\mu^\ast_{i,j} (\ell). \label{eq:asymptoticestimator}
\end{align}
Similarly, we have a.a.s.
\begin{align}
\hat{B}_{i,j} (\ell) \sim \frac{\omega}{n} \frac{   \int_{\mathcal{X}} h_{|\lambda_1| \epsilon}(d_r(\sigma_i,x) ) B_{x,\sigma_j} P(dx)  }{  \int_{\mathcal{X}} h_{|\lambda_1| \epsilon}(d_r(\sigma_i,x) ) P(dx)  }:=B^\ast_{i,j}.
\end{align}
The following proposition upper bounds the difference between the limit $\mu^\ast_{i,j}(\ell)$ (resp.\ $B^\ast_{i,j} (\ell)$) and $\mu_{\sigma_i,\sigma_j}(\ell)$ (resp.\ $B_{\sigma_i,\sigma_j}$).
\begin{proposition}\label{Prop2}
Suppose Assumption~\ref{assumptioncontinuous} holds. Then there exists a fraction of at least $(1-\sqrt{\epsilon_r})$ of all possible pairs  $(i,j)$ of nodes such
that
\begin{align}
B_{\sigma_i, \sigma_j} | \mu^\ast_{i,j} (\ell)-\mu_{\sigma_i, \sigma_j}(\ell) |  & \le  2 \psi(2|\lambda_1|\epsilon) + \frac{1}{|\lambda_1|^2\epsilon^2} \frac{ \sqrt{\epsilon_r}  }{ \int_{\mathcal{X} } h_{|\lambda_1| \epsilon}(d(\sigma_i,x) ) P(dx)}:= \eta, \; \forall \ell \in \calL, \nonumber \\
| B^\ast_{i,j}- \frac{\omega}{n} B_{\sigma_i, \sigma_j} | &  \le \frac{\omega}{n} \eta.
\end{align}
\end{proposition}
Applying Proposition~\ref{Prop2}, our theorem then follows.
\subsection{Proof of Theorem \ref{ThmNonReconstructionMultiple}}

Proof of Theorem \ref{ThmNonReconstructionMultiple} relies on a nice coupling between the local neighborhood of $\rho$ with a simple labeled Galton-Watson tree.
It is well-known that the local neighborhood of a node in the sparse graph is ``tree-like''. In the case with $r=2$, the coupling result is first studied in \cite{Mossel12} and generalized to the labeled tree in \cite{Jiaming13}. In this paper, we extend the coupling result to any finite $r \ge 2$.

Let $d=\omega \frac{a+(r-1)b}{r(a+b)}$ and consider a labeled Galton-Watson tree $\mathcal{T}$ with Poisson offspring distribution with mean $d$. The attribute of root $\rho$ is chosen uniformly at random from $\mathcal{X}$. For each child node, independently of everything else, it has the same attribute with its parent with probability $\frac{a}{a+(r-1)b}$ and one of $r-1$ different attributes with probability $\frac{b}{a+(r-1)b}$. Every edge between the child and its parent is independently labeled with distribution $\mu$ if they have the same attribute and with distribution $\nu$ otherwise.

The labeled Galton-Watson tree $\mathcal{T}$ can also be equivalently described as follows. Each edge is independently labeled at random according to the probability distribution $\mathbb{P}(\ell)=\frac{a \mu(\ell) + (r-1) b\nu(\ell)  }{ a + (r-1)b }$. The attribute of root $\rho$ is first chosen uniformly at random from $\mathcal{X}$. Then, for each child node, independently of everything else, it has the same attribute with its parent with probability $1-(r-1)\epsilon(\ell)$ and one of $r-1$ different attributes with probability $\epsilon(\ell)$, where
\begin{align}
\epsilon(\ell)=\frac{b \nu(\ell)}{a \mu(\ell)+ (r-1) b \nu(\ell) }. \label{eq:defepsilonell}
\end{align}

Recall that $G_R$ denote the neighborhood of $\rho$ in $G$ within distance $R$ and $\partial G_R$ denote the nodes at the boundary of $G_R$. Let $\mathcal{T}_R$ denote the tree $\mathcal{T}$ up to depth $R$ and $\partial \mathcal{T}_R$ denote the set of leaf nodes of $\mathcal{T}_R$. The following lemma similar to coupling lemmas in \cite{Mossel12} and \cite{Jiaming13} shows that $G_R$ can be coupled with the labeled Galton-Watson tree $\mathcal{T}_R$.
\begin{lemma} \label{LemmaCouplingTree}
Let $R=\theta \log n$ for some small enough constant $\theta>0$, then there exists a coupling such that a.a.s.\
$
(G_R,\sigma_{G_R})=(\mathcal{T}_R,\sigma_{\mathcal{T}_R}),
$
where $\sigma_{G_R}$ denote the node attributes on the subgraph $G_R$.
\end{lemma}

For the labeled Galton-Watson tree, we show that if $\omega< \omega_0$, then the attributes of leaf nodes are asymptotically independent with the attribute of root.
\begin{lemma} \label{ThmReconstructionTreeMultiple}
Consider a labeled Galton-Waltson tree $\mathcal{T}$ with $\omega< \omega_0$. Then as $R \to \infty$,
\begin{align}
 \forall x \in  \{1, \ldots, r\},\;\;\mathbb{P} ( \sigma_{\rho} =x | \mathcal{T}, \sigma_{\partial T_R} ) \to \frac{1}{r} \text{ a.a.s.} \nonumber
\end{align}
\end{lemma}
By exploiting Lemma \ref{LemmaCouplingTree} and Lemma \ref{ThmReconstructionTreeMultiple}, we give our proof of Theorem \ref{ThmNonReconstructionMultiple}.
By symmetry, $\mathbb{P}[\sigma_\rho=x |G, \sigma_v=y]=\mathbb{P}[\sigma_\rho=x' |G, \sigma_v=y]$ for $x, x' \neq y$ and $x \neq x'$.
Therefore, we only need to show that $\mathbb{P}[ \sigma_\rho=y | G,\sigma_v=y] \sim 1/r$ for any $y \in \mathcal{X}$ and
it further reduces to showing that
\begin{align}
\mathbb{P}[ \sigma_\rho=y | G,  \sigma_v=y,\sigma_{\partial G_R} ]  \sim 1/r.   \label{EqVarX}
\end{align}
Let $R=\theta \log n$ be as in Lemma \ref{LemmaCouplingTree} such that $G_R=o(\sqrt{n})$ and thus $v \notin G_R$ a.a.s.. Lemma 4.7 in \cite{Mossel12} shows that $\sigma_{\rho}$ is asymptotically independent with $\sigma_v$ conditional on $\sigma_{\partial G_R}$. Hence,
$
\mathbb{P}[ \sigma_\rho=y | G,  \sigma_v=y,\sigma_{\partial G_R} ]  \sim \mathbb{P} [\sigma_\rho= y |G,  \sigma_{\partial G_R} ].
$
Also, note that
$
\mathbb{P}(\sigma_\rho=y | G, \sigma_{\partial G_R} )= \mathbb{P}(\sigma_\rho=y | G_R, \sigma_{\partial G_R} ).
$
Lemma \ref{LemmaCouplingTree} implies that
$
\mathbb{P}(\sigma_\rho=y | G_R,\sigma_{G_R} ) \sim \mathbb{P}(\sigma_\rho=y | \mathcal{T}_R, \sigma_{\partial \mathcal{T}_R} ),
$
and by Lemma \ref{ThmReconstructionTreeMultiple}, $
 \mathbb{P}(\sigma_\rho=y | \mathcal{T}_R, \sigma_{\partial \mathcal{T}_R} ) \sim \frac{1}{r}.$
Therefore, equation (\ref{EqVarX}) holds.

\section*{Acknowledgment}
M.L. acknowledges the support of the French Agence Nationale de la Recherche (ANR) under reference ANR-11-JS02-005-01 (GAP project). J. X. acknowledges the support of NSF ECCS 10-28464.

\bibliographystyle{abbrv}
\bibliography{BibCommunityDetection}

\begin{appendix}

\section{Proof of Proposition \ref{PropSpectrum}  }
We first introduce notations used in the proof. Several norms on matrices will be used. The spectral norm of a matrix $X$ is denoted by $\|X\|$ and equals the largest singular value. The Frobenius norm of a matrix $X$ is denoted by $\|X\|_F$ and equals the square root of the sum of squared singular values. It follows that $\|X\|_F \le \sqrt{r} \|X\|_2$ if $X$ is of rank $r$. For vectors, the only norm that will be used is the usual $l_2$ norm, denoted as $\|x\|_2$.
Introduce a $n \times n$ matrix $\hat{A}$ defined by $\hat{A}_{ij}=K(\sigma_i,\sigma_j).$ Recall that $r$ is a fixed positive integer in Spectral Algorithm~\ref{alg:spectral}. Denote $r$ the largest eigenvalues of $\hat{A}$ sorted in decreasing order by $|\lambda'^{(n)}_1| \ge \cdots \ge |\lambda'^{(n)}_r|$ of $\hat{A}$. Let $v'_1,\ldots,v'_r\in\mathbb{R}^n$ be the corresponding eigenvectors with unit norm. An overview of the proof is shown in Fig.~\ref{fig:proofoutline}.
\begin{figure}[h!]
\begin{center}
\scalebox{1}{\begin{tikzpicture}[scale = 2, font = \small, thick]
\draw[<->] (-0.2, 0) node [left] {$\tilde{A}$ }-- (1, 0) node [right]{$\hat{A}$};
\draw[<->] (1.3,0)--(2.2,0) node [right] {$T$};
\node [above] at (0.4,0) {Lemma \ref{LemmaSpectral} and~\ref{LemmaDavisKahan} };
\node [above] at (1.8,0) {Lemma \ref{thm-koltchinskii}};
\node [below] at (-0.2,-0.07) {$\{ \lambda_k^{(n)}, v_k  \}_{k=1}^r$};
\node [below] at (1.2,-0.07) {$\{\lambda'^{(n)}_k,v'_k \}_{k=1}^r$};
\node [below] at (2.3,-0.10) {$\{ \lambda_k, \phi_k \}_{k=1}^r$};
\end{tikzpicture}}
\caption{Proof outline for showing the asymptotic correspondence between the spectrum of $\tilde{A}$ and that of $T$.}
\label{fig:proofoutline}
\end{center}
\end{figure}

Lemma \ref{thm-koltchinskii} follows from the results of~\cite{koltchinskii98} and their application as per Theorem 4 and Theorem 5 of \cite{vonluxburg-bousquet-belkin}.
\begin{lemma}\label{thm-koltchinskii}
Under our assumptions on operator $T$, for $k=1,2,\ldots, r+1$, almost surely $\frac{1}{n}\lambda'^{(n)}_k \sim \lambda_k$, and there exist choices of orthonormal eigenfunctions $\phi_k$ of operator $T$ associated with $\lambda_k$ such that $\lim_{n\to\infty}\sum_{i=1}^n(v'_k(i)-\frac{1}{\sqrt{n}}\phi_k(\sigma_i))^2=0$.
\end{lemma}
Lemma~\ref{LemmaSpectral} gives sharp controls for the spectral norm of random symmetric matrices with bounded entries initially developed by~\cite{Feige05} and extended by~\cite{mastom11} and~\cite{Chattergee12}.
\begin{lemma}\label{LemmaSpectral}
Let $M$ be a random symmetric matrix with entries $M_{ij}$ independent up to symmetry, $M_{ij}\in[0,1]$ and such that $\mathbb{E}[M_{ij}] =\omega/n$. If $\omega \ge  C \log n /n $ for a universal positive constant $C$, then for all $c>0$ there exists $c'>0$ such that with probability at least $1-n^{-c}$, one has
\begin{equation}
\| M-\mathbb{E}[M] \| \le c'\sqrt{\omega}.
\end{equation}
\end{lemma}
Lemma \ref{LemmaDavisKahan}, a consequence of the famous Davis-Kahan $\sin \theta$ Theorem \cite{Kahan70}, controls the perturbation of eigenvectors of perturbed matrices.
\begin{lemma}\label{LemmaDavisKahan}
For two symmetric matrices $M$, $M'$ and orthonormal eigenvectors $(u_1,\ldots,u_r)$ (respectively $u'_1,\ldots,u'_r)$ associated with the $r$ leading eigenvalues of $M$ (respectively $M'$), denoting $U=[u_1,\ldots,u_r]$, $U'=[u_1',\ldots u_r']$,  there exists an orthogonal $r \times r $ matrix $O$ such that
\begin{align}
\| U - U' O \| \le \frac{\sqrt{2} \| M-M'\| }{ |\theta_r| - |\theta_{r+1}| -\| M-M'\| }, \label{eq:DavisKahan}
\end{align}
where $\theta_k$ is the $k$-th largest eigenvalue of $M'$ in absolute value.
\end{lemma}

We omit proofs of lemmas which can be found in the mentioned literature. Next, we present the proof of our proposition. Applying Lemma~\ref{LemmaDavisKahan} to $M=\tilde{A}$ and $M'=(\omega/n)\hat{A}$, then we have $U=[v_1,\ldots, v_r]$, $U'=[v'_1,\ldots, v'_r]$, and $\theta_k=(\omega/n) \lambda'^{(n)}_k$ for $k=1,\ldots, r+1$.
By Lemma \ref{LemmaSpectral} and observing that $\mathbb{E}[\tilde{A}]=(\omega/n) \hat{A}$, it readily follows that $\| M -M'\| =O(\sqrt{\omega})$ with high probability. By Weyl's inequality, we have
$| \lambda_k^{(n)}- \theta_k| \le \|M-M'\|=O(\sqrt{\omega})$.  Moreover, by Lemma~\ref{thm-koltchinskii}, for $k=1, \ldots, r+1$, $\theta_k \sim \omega \lambda_k $. Hence, $\lambda_k^{(n)}/ \lambda_1^{(n)} = \lambda_k /\lambda_1 +O(1/\sqrt{\omega})$.
By assumption, $|\lambda_r| > |\lambda_{r+1}|$, and thus the right-hand side of~\eqref{eq:DavisKahan} is $O(1/\sqrt{\omega})$. Note that $U,U'O$ are of rank $r$, it follows that
\[
\| U - U' O \|_F \le \sqrt{2r} \| U - U' O \| = O(1/\sqrt{\omega}).
\]
Therefore, by Lemma~\ref{thm-koltchinskii}, there exist choices of orthonormal eigenfunctions $\phi_k$ of operator $T$ associated with $\lambda_k$ such that $\lim_{n\to\infty}\sum_{i=1}^n( v_k(i)-\frac{1}{\sqrt{n}}\phi_k(\sigma_i))^2=0$ for $k=1,\ldots, r$.

\section{Proof of Proposition~\ref{Prop2}}
The main idea of proof is to show that the pseudo-distance $d_r$ is close to distance $d$ in an appropriate sense. By definition, $d(x,x')\ge d_r(x,x')$ and moreover,
$$
\int_{\mathcal{X}^2}  [d^2(x,x')-d^2_r(x,x')]P(dx)P(dx') = \sum_{k>r} \lambda^2_k \int_{\mathcal{X}^2} \left(\phi_k(x)-\phi_k(x')\right)^2 \le  2 \sum_{k>r}\lambda_k^2=2 \epsilon_r.
$$
Define $d^2_{>r}(x,x')=d^2(x,x')-d^2_r(x,x')$. Markov's inequality entails that
$$
\int_{\mathcal{X}}P(dx')  \1{ \int_{\mathcal{X}} d^2_{>r}(x,x')P(dx) \ge 2 \sqrt{\epsilon_r}}\le \sqrt{\epsilon_r}.
$$


Note that the following inequalities hold
\begin{equation}\label{lip000}
0\le h_{\epsilon}(d_r(x,x'))-h_{\epsilon}(d(x,x')) \le \frac{1}{2\epsilon^2}\left[d^2(x,x')-d^2_r(x,x')\right]=\frac{d^2_{>r}(x,x')}{2\epsilon^2}.
\end{equation}
By the previous application of Markov's inequality, for a fraction of at least $1-\sqrt{\epsilon_r}$ of nodes $i$, it holds that
$$
\int_{\mathcal{X}} d^2_{>r}(x,\sigma_i)P(dx)\le 2\sqrt{\epsilon_r}.
$$
Combined with the previous Lipschitz property~(\ref{lip000}) and the definition of $ \mu^\ast$ given by~\eqref{eq:asymptoticestimator}, this entails that for a fraction of at least $1-\sqrt{\epsilon_r}$ nodes $i$, one has
$$
\frac{a}{b+\frac{\sqrt{\epsilon_r}}{|\lambda_1|^2\epsilon^2}} \le \mu^\ast_{i,j} (\ell) \le \frac{a+\frac{\sqrt{\epsilon_r}}{|\lambda_1|^2\epsilon^2} }{b},
$$
where we have introduced the following notations
\begin{align*}
a=\int_{\mathcal{X}} h_{|\lambda_1|\epsilon}(d(x,\sigma_i)\nu_{x,\sigma_j}(\ell)P(dx), \; b=\int_{\mathcal{X}}h_{|\lambda_1|\epsilon}(d(x,\sigma_i)B_{x,\sigma_j}P(dx).
\end{align*}
Define
\begin{align*}
a'= \int_{\mathcal{X} }h_{|\lambda_1| \epsilon}(d(\sigma_i,x) ) \nu_{\sigma_i,\sigma_j}(\ell) P(dx), \; b'=  \int_{\mathcal{X} } h_{|\lambda_1| \epsilon}(d(\sigma_i,x) )  B_{\sigma_i,\sigma_j}P(dx).
\end{align*}
Then, $\mu_{\sigma_i,\sigma_j}(\ell)=a'/b'$.
Note that for positive constants $c_1 \le c_2, c_3 \le c_4$, $|\frac{c_1}{c_2} - \frac{c_3}{c_4}| \le \frac{1}{c_4} (|c_1-c_3|+|c_2-c_4|)$. Hence,
$$
| \mu^\ast_{i,j} (\ell)- \mu_{\sigma_i,\sigma_j} (\ell)| \le \frac{ |a-a'| +  |b-b'| +\frac{\sqrt{\epsilon_r}}{|\lambda_1|^2\epsilon^2} }{ b' }.
$$
By assumption~\ref{assumptioncontinuous}, for all $x,x',y,y'$,
$$
|B_{x,y}-B_{x',y'}|\le \psi(d(x,x')+d(y,y'))
$$
and similarly for $\nu_{x,y}(\ell)$.
Therefore,
\begin{align*}
|a-a'| + |b-b'| &\le \int_{\mathcal{X} }h_{|\lambda_1| \epsilon}(d(\sigma_i,x) )\left[|\nu_{x,\sigma_i}(\ell)-\nu_{\sigma_i,\sigma_j}(\ell)|+|B_{x,\sigma_i}-B_{\sigma_j,\sigma_i}|\right]  P(dx) \\
& \le  2\psi(2|\lambda_1|\epsilon)\int_{\mathcal{X} }h_{|\lambda_1| \epsilon}(d(\sigma_i,x))P(dx).
\end{align*}
It follows that
\begin{align*}
B_{\sigma_i,\sigma_j} |\mu^\ast_{i,j} (\ell)- \mu_{\sigma_i,\sigma_j} (\ell)| & \le 2 \psi(2|\lambda_1|\epsilon) + \frac{\sqrt{\epsilon_r}}{|\lambda_1|^2\epsilon^2} \frac{ 1}{ \int_{\mathcal{X} } h_{|\lambda_1| \epsilon}(d(\sigma_i,x) ) P(dx)} =\eta.
\end{align*}
The right-hand side goes to zero as one lets successively $\epsilon_r$ then $\epsilon$ go to zero.
Similarly, we can show $ | B^\ast_{i,j}- \frac{\omega}{n} B_{\sigma_i, \sigma_j} |   \le \frac{\omega}{n} \eta.$

\section{Proof of Lemma~\ref{ThmReconstructionTreeMultiple}}
The proof technique is adapted from Section 5 in \cite{Mossel04}. Consider two distributions on the labeled Galton-Watson tree, one with the attribute of the root being $x$, and one with the attribute of the root being $y \neq x$. We can couple the two distributions in the following way: if the two distributions agree on the attribute of node $v$, then couple them together such that they also agree for all the children of $v$; if the two distributions do not agree on the attribute of node $v$, use the optimal coupling to make them agree as much as possible for each children of $v$. For each children $w$ with $L_{vw}=\ell$, it is easy to check that under the optimal coupling, the two distributions will not agree on the attribute of $w$ with probability $|1-r \epsilon(\ell) |$, where $\epsilon(\ell)$ is defined in~\eqref{eq:defepsilonell}. Hence, the non-coupled nodes grow as a branching process with the branching number given by
\[
\omega \sum_{\ell} \frac{a \mu(\ell) + (r-1) b \nu(\ell) }{r( a+b) }  |1-r \epsilon(\ell) | = \omega \tau.
\]
It is well known that if the branching number $\omega \tau <1$, the branching process will eventually die a.a.s. Thus as $ R \to \infty$, a.a.s.
\begin{align}
\mathbb{P} (  \sigma_{\partial T_R} | \mathcal{T},\sigma_{\rho} =x) =\mathbb{P} (  \sigma_{\partial T_R} | \mathcal{T},\sigma_{\rho} =y).  \nonumber
\end{align}
By Bayes' formula, the theorem follows.

\section{Bernstein Inequality}
\begin{theorem} \label{ThmBernstein}
Let $X_1, \ldots, X_n$ be independent random variables such that $| X_i | \le M$ almost surely. Let $\sigma_i^2=\text{Var}(X_i)$ and $\sigma^2= \sum_{i=1}^n \sigma_i^2$, then
\begin{align}
\mathbb{P} ( \sum_{i=1}^n X_i \ge t ) \le \exp \left (  \frac{-t^2}{ 2 \sigma^2 + \frac{2}{3} M t }\right). \nonumber
\end{align}
It follows then
\begin{align}
\mathbb{P} ( \sum_{i=1}^n X_i \ge \sqrt{2 \sigma^2 u} + \frac{2M u }{3}  ) \le e^{-u}. \nonumber
\end{align}
\end{theorem} 
\end{appendix}

\end{document}